\documentclass[11pt]{article}
\usepackage[latin1]{inputenc}
\usepackage{amsmath}
\usepackage{amsfonts}
\usepackage{amsthm}
\usepackage{amssymb,amsthm, stmaryrd, mathabx}
\usepackage{graphicx}
\usepackage{xcolor}
\usepackage{multicol}
\usepackage{soul}
\usepackage[normalem]{ulem}
\usepackage{float}
\usepackage[pagewise, displaymath, mathlines]{lineno}
\usepackage[displaymath, mathlines]{lineno}
\usepackage[numbers,sort&compress]{natbib}
\newcommand{\comment}[1]{}

\newtheorem{theorem}{Theorem}[section]

\newtheorem{proposition}[theorem]{Proposition}

\newtheorem{corollary}[theorem]{Corollary}
\newtheorem{definition}[theorem]{Definition}
\newtheorem{lemma}[theorem]{Lemma}

\newtheorem{example}[theorem]{Example}

\title{Recursive constructions of amoebas}

\begin{document}

\maketitle

\begin{center}

\begin{multicols}{2}

Adriana Hansberg\\[1ex]
{\small Instituto de Matem\'aticas\\
UNAM Juriquilla\\
Quer\'etaro, Mexico\\
ahansberg@im.unam.mx}\\[2ex]

\columnbreak

Amanda Montejano\\[1ex]
{\small UMDI, Facultad de Ciencias\\
UNAM Juriquilla\\
Quer\'etaro, Mexico\\
amandamontejano@ciencias.unam.mx}\\[4ex]

\end{multicols}

Yair Caro\\[1ex]
{\small Dept. of Mathematics\\
University of Haifa-Oranim\\
Tivon 36006, Israel\\
yacaro@kvgeva.org.il}\\

\end{center}

\vspace{2ex}

\begin{abstract} 
Global amoebas are a wide and rich family of graphs that emerged from the study of certain Ramsey-Tur\'an problems in $2$-colorings of the edges of the complete graph $K_n$ that deal with the appearance of unavoidable patterns once a certain amount of edges in each color is guaranteed. Indeed, it turns out that, as soon as such coloring constraints are satisfied and if $n$ is sufficiently large, then every global amoeba can be found embedded in $K_n$ such that it has half its edges in each color. Even more surprising, every bipartite global amoeba $G$ is unavoidable in every tonal-variation, meaning that, for any pair of integers $r, b$ such that $r + b $ is the number of edges of $G$, there is a subgraph of $K_n$ isomorphic to $G$ with $r$ edges in the first color and $b$ edges in the second. The feature that makes global amoebas work are one-by-one edge replacements that leave the structure of the graph invariant.  By means of a group theoretical approach, the dynamics of this feature can be modeled.  As a counterpart to the global amoebas that ``live'' inside a possibly large complete graph $K_n$, we also consider local amoebas which are spanning subgraphs of $K_n$ with the same feature. 

In an effort to highlight their richness and versatility, we present here three different recursive constructions of amoebas, two of them yielding interesting families per se and one of them offering a wide range of possibilities.\\

Keywords: Amoeba, edge-replacement, graph isomorphisms, interpolation, Ramsey-Tur\'an problems.
\end{abstract}

\section{Introduction}\label{intro}

Global amoebas first appeared in \cite{CHM3} (there they were just called ``amoebas'') where certain Ramsey-Tur\'an problems were considered, which dealt with the existence of a given graph with a prescribed color pattern in $2$-edge-colorings of the complete graph. More precisely, it is shown in \cite{CHM3} that, if $n$ is sufficiently large, and the edges of $K_n$ are colored each with one of two colors such that every color is sufficiently represented, then every global amoeba $G$ can be found embedded in $K_n$ such that it has half ($\pm 1$) its edges in each color. Graphs with this property are called \emph{balanceable}. Moreover, every bipartite global amoeba $G$ is \emph{omnitonal}, meaning that, provided that $n$ is sufficiently large and that the coloring of $K_n$ has enough edges in each color, then, for any pair of integers $r, b$ such that $r + b = e(G)$ (i.e. the number of edges of $G$), there is a subgraph of $K_n$ isomorphic to $G$ with $r$ edges in the first color and $b$ edges in the second. 

The feature that makes amoebas work are one-by-one replacements of edges, where, at each step, some edge is substituted by another such that an isomorphic copy of the graph is created. We call such edge substitutions \emph{feasible edge-replacements}.  Similar edge-operations have been studied, for instance, in \cite{CiZi,FrRi,FrHlRo,JaPaSc,Pun3,Ross}. 

As a counterpart to the global amoebas that ``live'' inside a possibly large complete graph $K_n$, we also consider local amoebas which are spanning subgraphs of $K_n$ with the same feature. To be more precise, consider, as introduced in \cite{CHLZ}, a family $\mathcal{F}$ of graphs, all of them having the same number of edges,  and a graph $H$.  $\mathcal{F}$  is called \emph{closed} in $H$ if,  for every two copies $F,  F'$ of members of $\mathcal{F}$ contained in $H$ as subgraphs,  there is a chain of graphs $H_1, H_2, \ldots, H_k$ in $H$ such that $F = H_1$,  $F' = H_k$, and, for $2 \le i \le  k$,  $H_i$ is a member of $\mathcal{F}$ and $H_i$ obtained from $H_{i-1}$ by a feasible edge-replacement.  Perhaps the most well-known closed family is the family of all spanning trees of a connected graph $H$, and the edge-replacement operation given above is in fact the basic operation in the exchange of bases in the cycle matroid $M(H)$ of $H$.  A graph $G$ is a \emph{global amoeba} precisely when $\{G\}$ is a closed family in $K_n$ (for $n$ large enough), and it is a \emph{local amoeba} if $\{G\}$ is a closed family in $K_{n(G)}$. Exactly this amoeba-property is the key role to the usefulness of amoebas in interpolation theorems in Graph Theory  and in zero-sum extremal problems \cite{CHLZ}, and in problems about forced patterns in $2$-colorings of the edges of $K_n$ \cite{CHM3}. For the interested reader, we refer to \cite{BaNiWh,HaMoPl,HaPl,LiLiTo,Pun1,Pun3,Zhou2} for more literature related to interpolation techniques in graphs.
 
In \cite{CHM-am}, a group theoretical framework is developed with which the dynamics of the feasible-edge replacements of graphs can be modeled. In order to present the results that are the aim of this paper, we will describe this framework formally but without going into detail. To understand the techniques more deeply, the reader is invited to consult \cite{CHM-am}. Two of the constructions that we expose in this paper are also matter from \cite{CHM-am}, where the full proofs can be read. Section 3.2, about the so-called ``double-rooted'' amoebas, is not part of \cite{CHM-am} and is fully contained and proven in the present manuscript. 

\section{Theoretical setting}\label{sec:theor_set}

Given a graph $G$,  we denote with $V(G)$, $E(G$), $n(G)$, and $e(G)$ its vertex set,  its edge set,  its order and its number of edges, respectively.  For integers $m$ and $n$ with $m<n$ we use the standard notation $[n]=\{1,2,...,n\}$ and $[m,n]=\{m,m+1,m+2,,...,n\}$. Let $S_n$ be the symmetric group, whose elements are permutations of $[n]$. The group of automorphisms of a graph $G$ is denoted by ${\rm Aut}(G)$. Thus, ${\rm Aut}(K_n)\cong S_n$ where $K_n$ is the complete graph of order $n$ and, for any graph $G$ of order $n$, ${\rm Aut}(G)\cong S$ for some $S\leq S_n$. Let $V = \{v_1, v_2, \ldots, v_n\}$ be the set of vertices of  $K_n$. Let $G$ be a spanning subgraph of $K_n$ defined by its edge set $E(G) \subseteq E(K_n)$ and let $L_G = \{ij \;|\; v_iv_j \in E(G)\}$, where we do not distinguish between $ij$ and $ji$. For each $\sigma \in S_n$, we define $\lambda_{\sigma} : V \to [n]$ as the labeling of the vertices of $K_n$ defined by $\lambda_{\sigma}(v_i) = \sigma(i)$ and consider 
the copy $G_{\sigma}$ of $G$ embedded in $K_n$ defined by 
\[E(G_{\sigma}) = \{v_{\sigma^{-1}(i)}v_{\sigma^{-1}(j)} \;|\; ij \in L_G\}.\]
Hence, each labeled  copy  of $G$ embedded in $K_n$ corresponds to a permutation  $\sigma \in S_n$ and vice versa. 
Observe that, for every (non-labeled) subgraph $G'$ of $K_n$ isomorphic to $G$ there are  $|{\rm Aut}(G)|$ different labelings of $V$ that correspond to $G'$, that is, the set  $\{\sigma \in S_n \,|\, G_{\sigma}=G'\}$ has  $|{\rm Aut}(G)|$ elements. Moreover,  $\{\sigma \in S_n \,|\, G_{\sigma}=G\}\cong {\rm Aut}(G)$. We will set 
\[A_G = \{\sigma \in S_n \,|\, G_{\sigma}=G\}.\]

It is important to note  that the set of labels of the edges of $G_{\sigma}$,  given by $L_{G_\sigma} = \{\sigma(i)\sigma(j) \;|\; v_iv_j \in E(G_{\sigma})\}$,  is the same for all $\sigma \in S_n$,  i.e. $L_{G_\sigma} = L_G$ for all $\sigma \in S_n$ (including $L_{G_{\rm id}} = L_G$). 
Moreover, the corresponding copies of the vertices and edges of $G$ in $G_{\sigma}$ are given by their labels: the copy of vertex $v_i$ of $G$ is the vertex of $G_{\sigma}$ having label $i$, while the copy of an edge $v_iv_j \in E(G)$ is the edge of $G_{\sigma}$ having label $ij$.

We denote by $\overline{G}$ the complement graph of $G$, that is, $V(\overline{G})=V(G)$ and $E(\overline{G})=\{uv \;|\; u, v \in V(G), uv \notin E(G) \}$.  Given $e\in E(G)$ and $e'\in E(\overline{G})$, we say that the graph $G-e+e'$ is obtained from $G$ by performing the \emph{edge-replacement} $e \to e'$. If $G-e+e'$ is a graph isomorphic to $G$, we say that the edge-replacement $e\to e'$ is \emph{feasible}.  We need also to consider the \emph{neutral edge-replacement} $\emptyset \to \emptyset$ as a feasible edge-replacement,  which is given when no edge is replaced at all.
Let $$R_G = \{rs \to kl \;|\;  G - v_rv_s + v_kv_l \cong G, rs \neq kl\} \cup\{\emptyset \to \emptyset\}$$ be the set of all feasible edge-replacements of $G$ given by their labels together with the neutral edge-replacement.  Let $R_G^* = R_G \setminus \{\emptyset \to \emptyset\}$.  Notice that,  since feasible edge-replacements different from the neutral one are defined by the labels of the edges,  any $rs \to kl \in R_G^*$ represents also a feasible edge-replacement of any copy $G_{\rho}$, $\rho \in S_n$.  Hence, clearly $R_{G_{\rho}} = R_G$ for any $\rho \in S_n$. 

Given a feasible edge-replacement, $rs \to kl\in R_G^*$, we will use the following notation   

\[S_G(rs \to kl)  = \{\sigma \in S_n \;|\; G_{\sigma} = G-v_rv_s+v_kv_l\},  \mbox{ and } S_G(\emptyset \to \emptyset)  = A_G.\]

We will use sometimes the notation $e \to e' \in R_G$ when we do not require to specify the indexes of the vertices involved in the edge-replacement, and it includes the possibility that $e \to e'$ is the neutral edge-replacement.

As it is developed in \cite{CHM-am}, it turns out that a feasible edge-replacement $e \to e' \in R_G$ in a copy $G_{\rho}$ of $G$ yields the copy of $G$ given by the permutation $\sigma 	\, \rho$, where we can choose any  $\sigma \in  S_G(e \to e')$. This way, we can define the group $S_G  := \left\langle \mathcal{E}_G\right\rangle$ generated by the permutations associated to all feasible edge-replacements, that is, by the set 
\[\mathcal{E}_G =  \bigcup_{e \to e'\in R_G}  S_G(e \to e'). \]

Clearly, $S_G$ acts on the set $\{ G_{\rho} \;|\; \rho \in S_n\} $ by means of $(\sigma,  G_{\rho}) 	\mapsto G_{\sigma \rho}$, where $\sigma \in S_G$ and $\rho \in S_n$. Observe that this action represents what happens when a series of feasible edge-replacements, represented by $\sigma$, is performed on a copy $G_{\rho}$ of $G$: the result is the graph $G_{\sigma	\rho}$.  

\begin{definition}\label{def:amoebas}
A graph $G$ of order $n$ is called a \emph{local amoeba} if $S_G = S_n$. That is, any labeled copy of $G$ embedded in $K_n$ can be reached, from $G$,  by a chain of feasible edge-replacements.
On the other hand, a graph $G$ is called \emph{global amoeba} if there is an integer $T \ge 0$ such that $G \cup tK_1$ is a local amoeba for all $t \ge T$.
\end{definition}

It is not difficult to convince oneself  that, for every $n\geq 2$, a path $P_n$ is both a local amoeba and  a global amoeba, while a cycle $C_n$ is neither a local amoeba nor a global amoeba, for any $n\geq3$.  In \cite{CHM-am}, several examples of amoebas are exhibited. It is important to note that there are examples in all possible combinations: graphs that are neither global nor local amoebas, graphs that are both local and global amoebas, and graphs that are just local or just global amoebas.

In \cite{CHM-am}, several structural properties of both graph families are presented and the very intrinsic relation between them is exhibited as well as their differences. Also some interesting constructions of both local and global amoebas are given (two of them are sketched here) and extremal global amoebas with respect to size, chromatic number and clique number are displayed. 

The following proposition is an easy consequence of the fact that, for any graph $G$, $S_G = S_{\overline{G}}$ (see  \cite{CHM-am}). 

\begin{proposition}\label{prop:basic}
A graph $G$ is a local amoeba if and only if $\overline{G}$ is a local amoeba. 
\end{proposition}

A major achievement in \cite{CHM-am} is the following characterization of global amoebas. It shows that, if we conceive a global amoeba as a graph $G$ embedded in $K_n$ that can be moved to any copy $G' \subseteq K_n$ of $G$ by means of a series of feasible edge-replacements, then $n$ can be as small as $n(G)+1$ (item (iii) of Theorem \ref{thm:eq}). This fact is surprising because one may intuitively think that, to be able to arrive to any copy of $G$, specially those with large intersection with the original graph $G$, one would need a lot of free room to move, which is not the case: we just need one vertex more than the order of $G$. Moreover, item (ii) of Theorem \ref{thm:eq} shows that, in order to detect if a graph $G$ of order $n$ is a global amoeba, we just need to show that every orbit of the action of $S_G$ on $[n(G)]$ contains an index $y$ belonging to a vertex of degree $1$.

\begin{theorem}\label{thm:eq}
Let $G$ be a non-empty graph  defined on the vertex set $V = \{v_1, v_2, \ldots, v_n\}$. The following statements are equivalent:
\begin{enumerate}
\item[(i)] $G$ is a global amoeba.
\item[(ii)] For each $x \in [n]$, there is a $y \in S_Gx$ such that $\deg_G(v_y) = 1$.
\item[(iii)] $G\cup K_1$ is a local amoeba.
\end{enumerate}
\end{theorem}

As a corollary of this theorem, we have the following facts.

\begin{corollary}\label{co:deg=0,1}
Let $G$ be a graph with minimum degree $\delta$.
\begin{enumerate}
\item[(i)] If $\delta\in \{0,1\}$ and $G$ is a local amoeba, then $G$ is a global amoeba.
\item[(ii)] If $\delta = 0$, then $G$ is a local amoeba if and only if $G$ is a global amoeba.
\end{enumerate}
\end{corollary}

\section{Recursive constructions of amoebas}

\subsection{The graph $H_n$}

In \cite{CHM-am}, we define the graph $H_n$ as follows. Let $V(H_n) = A \cup B$ such that, taking $q=\lfloor \frac{n}{2} \rfloor$, $A = \{v_1, v_2, \ldots, v_{q}\}$ and $B = \{v_{q+1}, v_{q+2}, \ldots, v_{q+{\lceil \frac{n}{2} \rceil}}\}$, where $B$ is a clique, $A$ is an independent set and adjacencies  between $A$ and $B$ are given by $v_i v_{q+j}\in E(H_n)$ if and only if $j \le i$, where $1\le i\le q$ and $1\le j\le \lceil \frac{n}{2} \rceil$. In the same paper, it is shown that $H_n$ can be defined recursively as we show in the next proposition.

\begin{proposition}
The graph $H_n$ can be given recursively by means of $H_1 = K_1$, and $H_n = \overline{H_{n-1} \cup K_1}$, for $n \ge 2$.
\end{proposition}

\begin{proof}
By definition, $H_1 = K_1$, and $H_2 \cong K_2$, being the latter the same as $\overline{H_1 \cup K_1}$. Now we will show that $H_n \cong \overline{H_{n-1} \cup K_1}$ for $n \ge 3$. Indeed, this comes from the fact that the set of all degree values in $H_{n-1} \cup K_1$ is $[n-2] \cup \{0\}$ yielding that the set of all degree values in $\overline{H_{n-1} \cup K_1}$ is $\{n-1 - d \; | 0 \le d \le n-2\} = [n-1]$. Hence, $H_n \cong \overline{H_{n-1} \cup K_1}$, for each $n \ge 2$. 
\end{proof}

By means of this recursion, and making use of the fact that $G$ is a local amoeba if and only if $\overline{G}$ is a local amoeba (Proposition \ref{prop:basic}) and that a local amoeba of minimum degree $1$ is also a global amoeba (Proposition \ref{co:deg=0,1}), it is not difficult to prove that $H_n$ is a local and global amoeba (see \cite{CHM-am} for details). The graph $H_n$ is interesting because it meets the upper bounds $e(G) \le \lfloor \frac{n^2}{4}\rfloor$, $\chi(G) \le \lfloor \frac{n}{2}\rfloor+1$, and $\omega(G) \le \lfloor \frac{n}{2}\rfloor+1$ for global amoebas $G$ of order $n$ and minimum degree $\delta(G) = 1$. That is, $H_n$ is a global amoeba that is as dense as possible and has the highest possible clique and chromatic numbers. We gather these facts in the following proposition, whose proof can be found in \cite{CHM-am}.

\begin{proposition}\label{prop:Hn_GA&LA}
For every $n \ge 2$, $H_n$ is a global and local amoeba with $\delta(G)=1$, $e(H_n) = \lfloor \frac{n^2}{4}\rfloor$ and  $\omega(H_n) = \lfloor \frac{n}{2}\rfloor+1$.
  \end{proposition}

\subsection{Double-rooted amoebas}

In this section, we will introduce a way of gluing global amoebas such that the result is again a global amoeba.\\

The following two lemmas are quite technical but are the key elements for constructing amoebas from smaller ones as we will do further on. 

\begin{lemma}\label{la:subgraphs}
Let $G$ be a graph on vertex set $V = \{v_i \;|\; i \in [n]\}$. Let $G = G' \cup G''$ for two subgraphs $G'$ and $G''$ with respective vertex sets $V'$ and $V''$. Let $J'$ and $J''$ be the sets of indexes of the vertices in $V'$ and $V''$, respectively, and let $I = J' \cap J''$. If there is a $\sigma \in \mathcal{E}_{G'} \cap \bigcap_{j \in I}{\rm Stab}_{S_{G'}}(j)$, then the permutation
\[
\widehat{\sigma}(i) = \left\{ \begin{array}{ll} 
\sigma(i),& \mbox{for } i \in J'\setminus I \\
i,& \mbox{for } i \in J''
\end{array}
\right.
\]
is in $\mathcal{E}_{G}$.
\end{lemma}

\begin{definition}\label{def:composition}
Let $G$ be a graph on vertex set $V = \{v_i \;|\; i \in [n]\}$ and $H$ another graph of order $m$ provided with a root. Let $I = \{i_1, i_2, \ldots, i_k\}\subseteq [n]$, $N = n + k (m-1)$ and let $[N] = [n] \cup \bigcup_{\ell = 1}^k J_{i_\ell}$ be a partition of $[N]$ such that $|J_{i_\ell}| = m-1$ for all $1 \le \ell \le k$. We define $G *_I H$ as the graph consisting of $G$ and copies $H_{i_1}, H_{i_2}, \ldots, H_{i_k}$ of $H$ such that $V(H_{i_\ell}) = \{v_i \;|\; i \in J_{i_\ell} \cup\{i_{\ell}\}\}$, for $1 \le \ell \le k$. If $I = [n]$, we set $G * H$ instead of $G *_{[n]} H$. Finally, we define $H^1 = G$ and, for $k \ge 2$, let $H^k := H^{k-1}*H$ as the \emph{$k$-th power graph} of $H$.
\end{definition}

\begin{lemma}\label{la:expansion}
Let $G$ be a graph on vertex set $V = \{v_i \;|\; i \in [n]\}$ and $H$ another graph of order $m$ provided with a root. Let $I = \{i_1, i_2, \ldots, i_k\}\subseteq [n]$, $N = n + k (m-1)$ and let $[N] = [n] \cup \bigcup_{\ell = 1}^k J_{i_\ell}$ be a partition of $[N]$ such that $|J_{i_\ell}| = m-1$ for all $1 \le \ell \le k$.  For any $\ell, \ell' \in [k]$, $\ell \neq \ell'$, let $\varphi_{i_\ell, i_{\ell'}} : J_{i_\ell} \to J_{i_{\ell'}}$ be the bijection given by an isomorphism between $H_{i_\ell}$ and $H_{i_{\ell'}}$ that sends $v_{i_{\ell}}$ to $v_{i_{\ell'}}$. If $\sigma \in \mathcal{E}_G \cap {\rm Stab}_{S_G}(I)$, then the permutation 
\[
\widetilde{\sigma}(i) = \left\{ \begin{array}{ll} 
\sigma(i),& \mbox{for } i \in [n] \\
\varphi_{i_\ell, \sigma(i_{\ell})}(i),& \mbox{for } i \in J_{i_\ell}, \ell \in [k] 
\end{array}
\right.
\]
is in $\mathcal{E}_{G *_I H}$.
\end{lemma}

To proceed with the last construction, we will need some terminology.

\begin{definition}
Let $G$ be a graph on $n$ vertices $v_1, v_2, \ldots, v_n$ rooted on $v_k$, for some $k \in [n]$. We call $G$ \emph{stem-transitive} if there is a set $S \in {\rm Stab}_{S_G}(k)$ that acts transitively on $[n] \setminus \{k\}$. Moreover, a vertex $v_j$, with $j \neq k$, such that there is a permutation $\varphi \in A_G$ with $\varphi(k) = j$ is called a \emph{root-similar vertex}.
\end{definition}

Observe that, a rooted graph $G$ with $\delta(G) = 1$ such that it is stem-transitive and has a root-similar vertex is a global amoeba by means of Theorem \ref{thm:eq}. We can then talk about stem-transitive (global or local) amoebas with a root-similar vertex. We recall here that a local amoeba with minimum degree $1$ is also a global amoeba (see Corollary \ref{co:deg=0,1}).

\begin{definition}
A rooted graph $G$ with $\delta(G) = 1$ such that it is stem-transitive and has a root-similar vertex is called a \emph{double-rooted global amoeba}. If, additionally, $G$ is a local amoeba, then $G$ is called a \emph{double-rooted local amoeba}.
\end{definition}

An example of such double-rooted local amoebas are the paths, where the root can be any vertex with the exception of the central vertex in odd-order paths. In \cite{CHM-am}, it was proved that paths are local amoebas.

\begin{proposition}\label{prop:paths}
A path $P = v_1v_2\ldots v_n$ on $n \ge 2$ vertices with root $v_k$ is a double-rooted local amoeba if and only if $n$ is even or if $n$ is odd and $k \neq \frac{n+1}{2}$. 
\end{proposition}

\begin{proof}
We will prove first that, for any $j \in [n-1]$, there is a set $S \subset \bigcap_{j+1 \le i \le n}{\rm Stab}_{S_P}(i)$ that acts transitively on $[j]$. We proceed by induction on $j$. If $j = 1$, there is nothing to prove. If $j = 2$ (and thus $n \ge 3$), consider the permutation  $(1\; 2) \in S_P( 2 \; 3 \to 1 \; 3)$ associated to the feasible edge-replacement  $2 \; 3 \to 1 \; 3$, showing that the set $\{(1\; 2)\}$, which is contained in $\bigcap_{3 \le i \le n}{\rm Stab}_{S_P}(i)$, acts transitively on $[2]$. For $n \le 3$, we have already finished. So we can assume that $n \ge 4$. Suppose now that $j \in [2, n-1]$ and that there is a set $S \subset \bigcap_{j+1 \le i \le n}{\rm Stab}_{S_P}(i)$ that acts transitively on $[j]$. Take the feasible edge-replacement $j+1 \; j+2 \to 1 \; j+2$ and observe that the permutation $\sigma$ defined by $\sigma(i) = j+2-i$ for $i \in [j+1]$ and $\sigma(i) = i$ for $i \in [j+2, n]$ belongs to $S_P(j+1 \; j+2 \to 1 \; j+2)$ and also to $\bigcap_{j+2 \le i \le n}{\rm Stab}_{S_P}(i)$. Hence, $S' = S \cup \{\sigma\}$ is a set of permutations in $S_P$ that acts transitively on $[j+1]$. By symmetry, we can conclude also that, for any $j \in [2,n]$ there is a set $S \subset \bigcap_{1 \le i \le j-1}{\rm Stab}_{S_P}(i)$ that acts transitively on $[j, n]$.

Now we turn to prove the statement of this proposition. Clearly, $v_k$ has a root-similar vertex if and only if $n$ is even or if $n$ is odd and $k \neq \frac{n+1}{2}$. Hence, the only if part of the proposition is clear. For the converse, we only need to prove that $P$ rooted on $v_k$ is stem-transitive. If $k = 1$ or $k = n$ we are done with what was proved above. So we assume that $k \in [2, n-1]$. Let $S \subset \bigcap_{k+1 \le i \le n}{\rm Stab}_{S_P}(i)$ and $S'\bigcap_{1 \le i \le k-1}{\rm Stab}_{S_P}(i) $ be sets that act transitively on $[k-1]$ and $[k+1, n]$, respectively. By symmetry, we can assume, without loss of generality, that $k \le \frac{n}{2}$. Consider now the feasible edge-replacement $(2k-1 \; 2k \to 1 \; 2k)$ with associated permutation $\tau \in S_P(2k-1 \; 2k \to 1 \; 2k)$ defined by $\tau(i) = 2k-1$ for $i \in [2k-1]$, and $\tau(i) = i$ for $i \in [2k, n]$. Observe that $\tau \in {\rm Stab}_{S_P}(k)$.  Moreover, since $S$ acts transitively on $[k-1]$ and $S'$ does the same on $[k+1, n]$, and $\tau(1) = 2k-1$ with $1 \in [k-1]$ and $2k-1 \in [k+1, n]$, it follows that $S \cup S' \cup \{\tau\}$ is a set of permutations acting transitively on $[n] \setminus \{k\}$, implying that $P$, rooted on $v_k$, is stem-transitive. 

Hence, as paths are local amoebas, we have proved that $P$, with root $v_k$, is a double-rooted local amoeba if and only if $n$ is even or if $n$ is odd and $k \neq \frac{n+1}{2}$. 
\end{proof}

\begin{proposition}\label{prop:Hn}
The graph $H_n$ rooted on a vertex of degree $\lfloor \frac{n}{2}\rfloor$ is a double-rooted global and local amoeba.
\end{proposition}

\begin{figure}[h]
\begin{center}
	\includegraphics[scale=0.6]{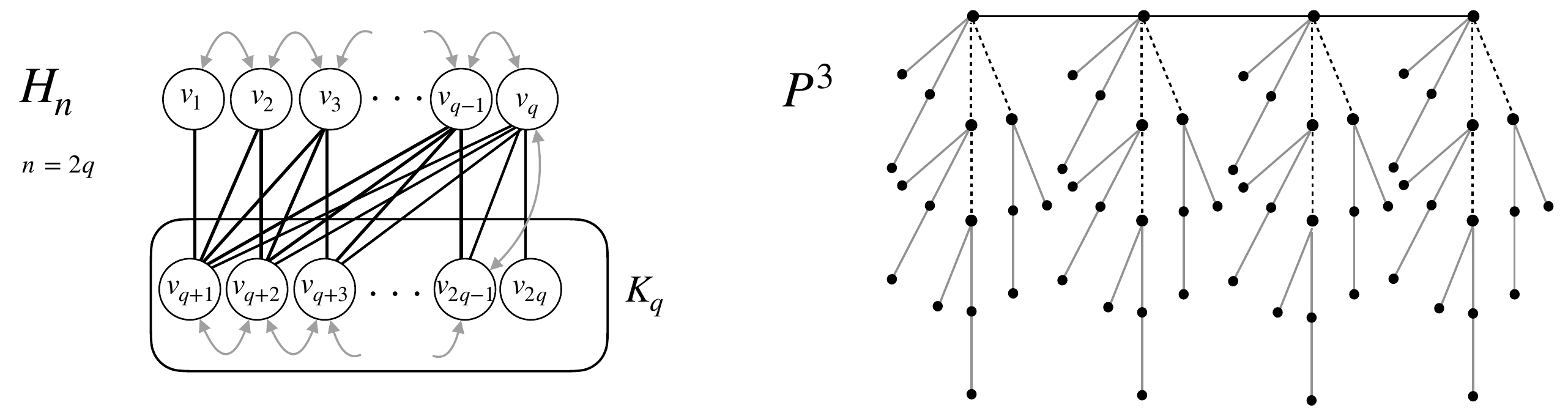}
	\caption{\small Sketch of the proof of Proposition \ref{prop:Hn}, and the graph $P^3$, where $P$ is a $P_4$ rooted on one of the vertices of degree $2$.}
	\label{fig:H_n&P^3}
\end{center}	
\end{figure}

\begin{proof}
We will prove only the case that $n$ is even as the odd case is completely analogous. Let $V(H_n) = A \cup B$ such that, taking $q=\frac{n}{2} $, $A = \{v_1, v_2, \ldots, v_{q}\}$ and $B = \{v_{q+1}, v_{q+2}, \ldots, v_{2q} \}$, where $B$ is a clique, $A$ is an independent set, and adjacencies between $A$ and $B$ are given by $v_i v_{q+j}\in E(H_n)$ if and only if $j \le i$, where $1\le i, j\le q$. Clearly, both $v_q$ and $v_{2q}$ are the only two vertices having degree $q=\frac{n}{2}$ and they are similar to each other. We designate $v_{2q}$ as the root of $H_n$. Now we will prove that the so rooted $H_n$ is stem-transitive. Indeed, the feasible edge replacements $( v_{i+1}\; v_{q+i+1} \to v_i \;v_{q+i+1}\;)$, for $i \in [q-1]$, and $( v_{i-q}\; v_i \to v_{i-q} \; v_{i+1}\;)$, for $i \in [q+1, 2q-2]$, lead to the associated permutations $(i \; i+1)$, for $i \in [2q-2] \setminus \{q\}$, which, together with the permutation $(q \; 2q-1)$ that is obtained via the feasible edge replacement $( v_{q-1}\; v_{2q-1} \to v_{q-1} \;v_q\;)$, shape a set of permutations that act transitively on $[n] \setminus \{2q\}$ (see Figure \ref{fig:H_n&P^3} for an illustration of this action). Hence, $H_n$ is stem-transitive. We can conclude now that $H_n$ is a double-rooted global and local amoeba.
\end{proof}

The following theorem provides a way of constructing global amoebas by means of smaller global amoebas, being one of the ingredients a double-rooted  global amoeba.

\begin{theorem}\label{thm:G*H}
Let $G$ be a non-empty global amoeba, and let $H$ be a double-rooted global amoeba. Then $G * H$ is a global amoeba.
\end{theorem}

\begin{proof}
Let $n = n(G)$ and $m = n(H)$. Let $N = nm$, which is the order of $G * H$. We set $[N] = [n] \cup \bigcup_{j = 1}^n J_j$ as the partition of $[N]$, where $V(H_j) = \{v_i \;|\; i \in J_j \cup\{j\}\}$, for $1 \le i \le n$, and $V(G) = \{v_i \;|\; i \in [n]\}$.
In view of Theorem \ref{thm:eq}, we will show that, for all $j \in [N]$, there is a permutation $\sigma \in S_{G * H}$ such that  $\deg_{G*H}(v_{\sigma(j)}) = 1$. We distinguish into three different cases.\\

\noindent
{\it Case 1. Suppose that $j \in J_k$ for some $k \in [n]$.}\\
Since $\delta(H_k) = 1$ and $H_k$ is stem-transitive, there is a permutation $\sigma \in S_{H_k}$ such that $\deg_{H_k}(v_{\sigma(j)}) = 1$ with $\sigma(j) \neq k$ (observe that if $v_k$ has degree~$1$ in $H_k$?, then its root-similar vertex, which has index in $J_k$, has also degree~$1$). Consider now $\widehat{\sigma} \in S_{H_k}$ such that
\[
\widehat{\sigma}(i) = \left\{ \begin{array}{ll} 
\sigma(i),& \mbox{for } i \in J_k, \\
i,& \mbox{for } i \in [N] \setminus J_k,
\end{array}
\right.
\]
which, by Lemma \ref{la:subgraphs}, belongs to $\mathcal{E}_{G*H}$. Then we have clearly that $\deg_{G*H}(v_{\widehat{\sigma}(j)}) = \deg_{H_k}(v_{\sigma(j)}) = 1$ and we are done. \\

\noindent
{\it Case 2. Suppose that $j \in [n]$ such that $\deg_G(v_j) = 1$.}\\
Let $N_G(v_j) = \{v_k \}$. Let $\varphi \in A_{H_j}$ such that, with $\ell = \varphi^{-1}(j)$, $v_{\ell}$ is the root-similar vertex of $v_j$. Then we know that $\ell \in J_j$. Consider the feasible edge-replacement $k \; j \to k \; \ell$ in $G*H$ with associated permutation $\tau \in S_{G*H}(k \; j \to k \; \ell)$ defined by
\[
\tau(i) = \left\{ \begin{array}{ll} 
\varphi(i),& \mbox{for } i \in J_j \cup \{j\}, \\
i,& \mbox{else}.
\end{array}
\right.
\]
Since $\tau(j) \in J_j$, by Case 1 there is a permutation $\sigma \in S_{G*H}$ such that $\deg_{G*H}(v_{\sigma\tau(j)}) = 1$. Hence, we have finished in this case, too.\\

\noindent
{\it Case 3. Suppose that $j \in [n]$ such that $\deg_G(v_j) \neq 1$.}\\
Since $G$ is a global amoeba, by Theorem \ref{thm:eq} there is a $\rho \in S_G$ such that $\deg_G(v_{\rho(j)}) = 1$. By previous case, there is a permutation $\sigma \in S_{G*H}$ such that $\deg_{G*H} (v_{\sigma \rho(j)}) = 1$, and we are done.\\

Since we have checked all possible cases, the proof is complete.
\end{proof}

\begin{proposition}
If $G$ and $H$ are both double-rooted global amoebas, then $G * H$ is again a double-rooted global amoeba with the same root as $G$.
\end{proposition}

\begin{proof}
Let $V(G) = \{v_1, v_2, \ldots, v_n\}$ and suppose that $v_k$ is the root of $G$ and $v_{k'}$ its root-similar vertex. Let $\psi \in A_G$ such  that $\psi(v_k) = v_{k'}$. Then $\psi$ can be extended to an element $\widetilde{\psi} \in A_{G * H}$ by means of Lemma \ref{la:expansion} (taking $I = [n]$ and noting that ${\rm Stab}_{S_G}([n]) = S_G$.). Hence, if we take $v_k$ as the root of $G * H$, then $\widetilde{\psi} (v_k) = \psi(v_k) = v_{k'}$ is a root-similar vertex to $v_k$. We now need to show that $G * H$ rooted on $v_k$ is stem-transitive. To this aim, we use the fact that $G$ is itself stem-transitive and that it contains a vertex different from the root that has degree $1$ and, proceeding as in the proof of Theorem \ref{thm:G*H}, we can conclude that $G*H$ is stem-transitive. Hence, $G * H$ is again a double-rooted amoeba with the same root as $G$.
\end{proof}

If $G$ is a double-rooted amoeba, we can apply Theorem \ref{thm:G*H} to create the double-rooted amoeba $G*G$ and repeat the process indefinitely.

\begin{corollary}
Let $G$ be a double-rooted global amoeba with root $v$.  Then $G^k$ is a double-rooted global amoeba with root $v$.
\end{corollary}

\begin{example}
Consider any non-empty global amoeba $G$ on $n$ vertices, where $V(G) = \{v_1,v_2,\ldots, v_n\}$. Let $P$ a path on $m$ vertices with rooted on some vertex except for the central vertex if $n$ is odd. Let $H$ be isomorphic to the graph $H_m$ and rooted on a vertex of degree $\lfloor \frac{m}{2}\rfloor$. By Propositions \ref{prop:paths} and \ref{prop:Hn}, $P$ and $H$ are double-rooted local amoebas. Hence, Theorem \ref{thm:G*H} yields that $G * P$, $G *H$ are global amoebas, while $P^k$ and $H^k$ are both double-rooted global amoebas, where $k \ge 1$ (see Figure \ref{fig:H_n&P^3} for an illustration of $P^3$, where we take one of the vertices of degree $2$ as the root). 
\end{example}





\subsection{Fibonacci amoeba-trees}

\begin{figure}[h]
\begin{center}
	\includegraphics[scale=0.3]{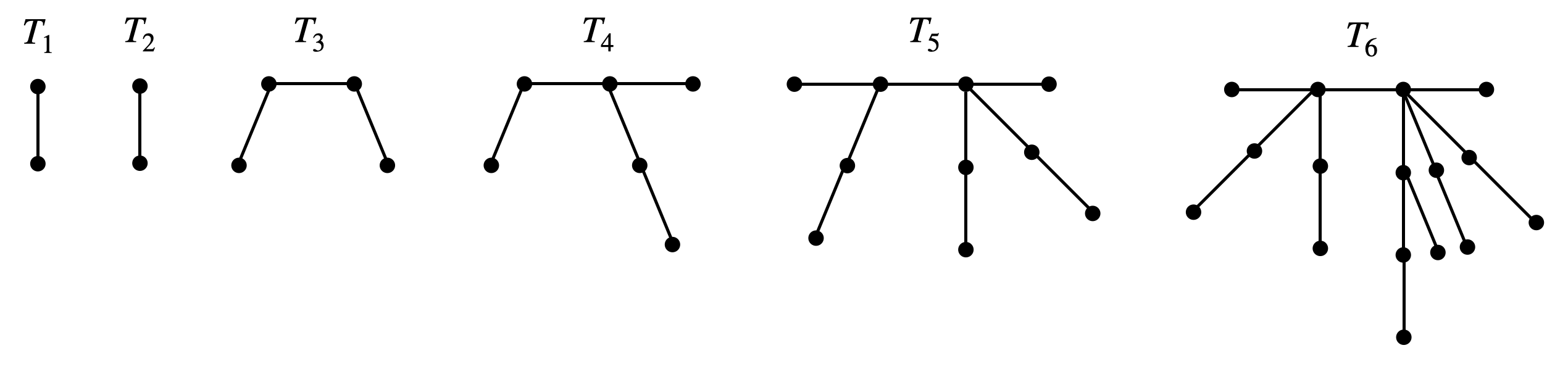}
	\caption{\small Fibonacci amoeba-trees $T_i$, $1 \le i \le 6$.}
	\label{fig:Fib_amoebas}
\end{center}	
\end{figure}

Next we will define a family of trees $(T_i)_{i \ge 1}$ that is constructed by means of a Fibonacci recursion. Let $T_1 = T_2 = K_2$. For $i \ge 2$, we define $T_{i+1}$ as the tree consisting of one copy $T$ of $T_{i-1}$ and one copy $T'$ of $T_i$, where a vertex of maximum degree of $T$ is joined to a vertex of maximum degree of $T'$ by means of a new edge (see Figure \ref{fig:Fib_amoebas}). Observe that $n(T_i) = 2F_i$, being $F_i$ the $i$-th Fibonacci number. Note also that, for $i \ge 4$, $T_i$ has only one vertex of maximum degree, which we will be the root of $T_i$. For the case that $i \le 3$, we will designate one of the vertices of maximum degree as the root of $T_i$.

\begin{theorem}\label{thm:Ti_GA}
$T_i$ is a global amoeba for all $i \ge 1$.
\end{theorem}

\begin{proof}{(Sketch)}
Let $T$ be a tree isomorphic to $T_i$. Let $J$ be the set of indexes of the vertices of $T$, i.e. $V(T) = \{v_k \;|\; k \in J\}$ and let $c \in J$ such that $v_c$ has maximum degree in $T$. We designate $v_c$ as the root of $T$. The proof consists of two parts: the first, where we show  that $T_i$ is stem-transitive (Claim 1); the second, where we prove that there is a feasible edge-replacement $e \to e'$ with an associated permutation in $S_T(e \to e')$ that does not fix $c$ (Claim 2). With these two facts, one can conclude that $S_T$ acts transitively on $J$, which in turn implies that, for every vertex $v_k$, there is a $\sigma \in S_T$ such that $\deg_T(v_{\sigma(k)}) = 1$. Hence, by Theorem \ref{thm:eq} it follows that $T$ is a global amoeba. \\

\noindent
{\it Claim 1. For all $i \ge 1$, $T_i$ is stem-transitive.}\\
We proceed by induction on $i$. If $i \le 3$,  $T$ is a path and the conclusion that $T$ is stem-transitive follows from Proposition \ref{prop:paths}. If $i = 4$, then let $T$ be the tree built from the path $v_4v_3v_1v_2 \cong T_3$ and a $T_2 \cong K_2$, given by $v_5v_6$, and the edge $v_1v_5$ joining both trees. Clearly, the only maximum degree vertex is $v_1$ and thus $c = 1$. Then the feasible edge-replacements $34 \to 24$ and $13 \to 14$ give respectively the permutations $(2 \, 3)$ and $(3 \, 4)$, which together with the automorphism  $(3 \, 5) (4 \, 6)$, act transitively on $[5] \setminus \{1\} = J \setminus \{c\}$ leaving $c = 1$ fixed. Hence,  $T_4$ is stem-transitive.

Now suppose that $i \ge 4$ and that we have proved the above statement for integer values at most $i$. Let $T \cong T_{i+1}$. For a subset $X \subset J$, we define $V_X = \{v_x \;|\; x \in X\}$ and $T_X = T[V_X]$. Let $J = U \cup W$ be a partition of $J$ such that $T_U \cong T_{i-1}$ and $T_W \cong T_i$. Further, let $U = A \cup B$ and $W = C \cup D$ be partitions such that $T_A \cong T_{i-3}$, $T_B \cong T_{i-2}$, $T_C \cong T_{i-1}$, and $T_D \cong T_{i-2}$. By construction, $v_c$ is the root of $T_C$. Let $a, b, d \in J$ be such that $v_a, v_b, v_d$ are the roots of $T_A, T_B$, and $T_D$, respectively.  Notice that $v_av_bv_cv_d$ is a path of length $4$ in $T$. See Figure \ref{fig:TreeT} for a sketch.

\begin{figure}[h]
\begin{center}
	\includegraphics[scale=0.4]{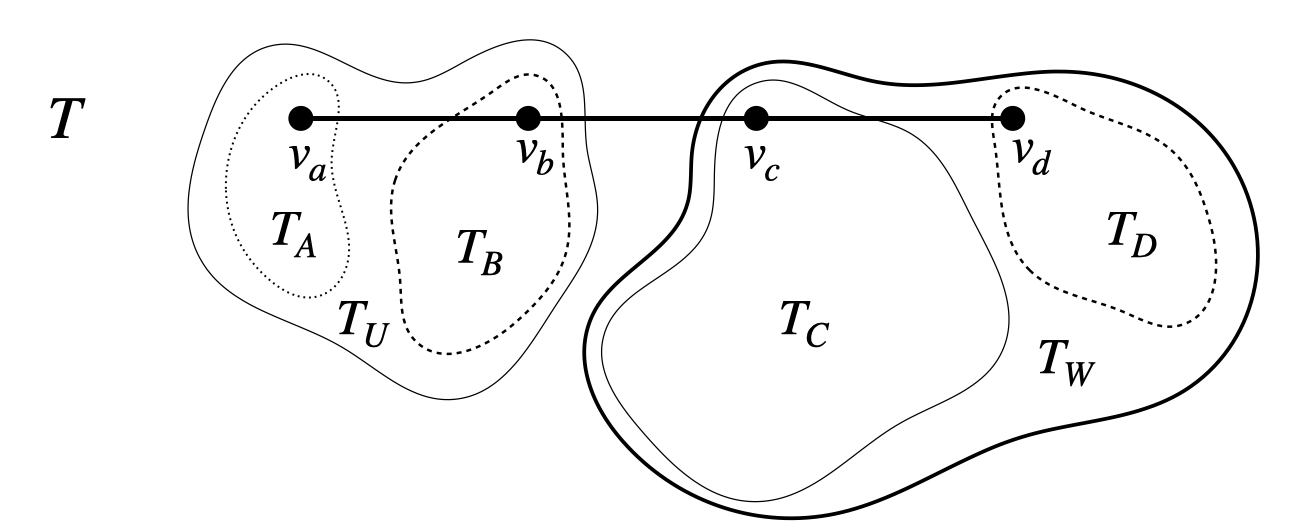}
	\caption{\small Sketch of the tree $T \cong T_{i+1}$ with its subtrees $T_U \cong T_{i-1}$ and $T_W \cong T_i$, and subsubtrees $T_A \cong T_{i-3}$, $T_B \cong T_{i-2}$, $T_C \cong T_{i-1}$, and $T_D \cong T_{i-2}$.}
	\label{fig:TreeT}
\end{center}	
\end{figure}

By the induction hypothesis, $T_U$ and $T_W$ are stem-transitive, where $v_b$ is the root of $T_U$, and $v_c$ the root of $T_W$. By means of Lemma \ref{la:subgraphs}, the so obtained permutations sets $Q_U \subseteq S_{T_U}$ and $Q_W \subseteq S_{T_W}$ that act transitively on $U \setminus \{ b\}$ and, respectively, on $W \setminus \{c\}$ can be used to produce permutations sets $\widehat{Q}_U, \widehat{Q}_W \subseteq S_T$ that inherit the transitive action on the corresponding set. 

Consider now the tree $T(B,D)$ that is obtained by identifying all vertices from $V_B$ with vertex $v_b$ and all vertices from $V_D$ with vertex $v_d$, i. e. we contract the sets $V_B$ and $V_D$ each to a single vertex. Observe  that $ab \to ad$ is a feasible edge-replacement in $T(B,D)$ with $\tau = (b \, d) \in S_{T(B,D)}(ab \to ad)$, and that $T \cong T(B,D) *_{\{b,d\}} T_{i-2}$. Considering that $T_B \cong T_D \cong T_{i-2}$, we can use now Lemma \ref{la:expansion} to obtain the corresponding permutation $\widetilde{\tau}  \in S_T(ab \to ad)$ which leaves $c$ fixed. It follows now easily that
\[S = \widehat{Q}_U \cup \widehat{Q}_W \cup \{ \widetilde{\tau} \}\]
is such that $\langle S \rangle$ acts transitively on $J \setminus \{c\}$. \\

\noindent
{\it Claim 2. There is a permutation $\widetilde{\rho} \in S_T$ such that $\langle S \cup \{\widetilde{\rho}\} \rangle$ acts transitively on $J$.}\\
Since we know already that $\langle S \rangle$ acts transitively on $J \setminus \{c\}$, we just need to find a $\widetilde{\rho}\in S_T$ with $\widetilde{\rho}(c) \neq c$. Indeed, there is such a permutation $\widetilde{\rho}$, namely one produced by the feasible edge-replacement $cd \to bd$ in $T$, which, by Lemma \ref{la:expansion}, can be obtained by means of the permutation $(b \, c) \in S_{T(U,C)}(cd \to bd)$.\\
\end{proof}

\section*{Acknowledgments}
We would like to thank BIRS-CMO for hosting the workshop Zero-Sum Ramsey Theory: Graphs, Sequences and More 19w5132, in which the authors of this paper were organizers and participants, and where many fruitful discussions arose that contributed to a better understanding of these topics.

\end{document}